\documentclass[12pt]{article}

\usepackage{amsmath,amssymb}		
\usepackage{float}

\usepackage{indentfirst,misccorr}	

\newcommand{\R}{\mathbb{R}}
\newcommand{\Q}{\mathbb{Q}}

\newtheorem{theorem}{Theorem}
\newtheorem{lemma}{Lemma}
\newtheorem{proposition}{Proposition}

\newtheorem{definition}{Definition}

\renewcommand{\le}{\leqslant}
\renewcommand{\ge}{\geqslant}
\renewcommand{\leq}{\leqslant}
\renewcommand{\geq}{\geqslant}

\title{On the chromatic numbers of small-dimensional Euclidean spaces}
\date{}
\author{Danila Cherkashin\footnote{Saint Petersburg State University, Faculty of Mathematics and Mechanics; Moscow Institute of Physics and Technology, 
Laboratory of Advanced Combinatorics and Network Applications.}, Anatoly Kulikov\footnote{Saint Petersburg State University, Faculty of Physics;
The University of Queensland, School of Mathematics and Physics.}, 
Andrei Raigorodskii\footnote{Moscow Institute of Physics and Technology, Faculty of Innovations and High Technology, Department of
Discrete Mathematics and Laboratory of Advanced Combinatorics and Network Applications; Moscow State University, Mechanics and Mathematics Faculty, 
Department of Mathematical Statistics 
and Random Processes; Buryat State University, Institute of Mathematics and Informatics.}}

\begin{document}

\maketitle

\begin{abstract}

This paper is devoted to the study of the graph sequence $G_n = (V_n, E_n)$, where $V_n$ is the set of all vectors $v \in \mathbb {R}^n$ with 
coordinates in $\{-1,0,1\}$ such that $|v| = \sqrt{3}$ and $E_n$ consists of all pairs of vertices with scalar product $1$. We find the exact value of the 
independence number of $G_n$. As a corollary we get new lower bounds on $\chi (\mathbb{R}^n)$ and $\chi (\mathbb{Q}^n)$ for small values of $n$.

\end{abstract}

\section{Introduction}                                                                                

Let $ {\mathbb R}^n $ be the standard Euclidean space, where the distance between any two points $ x, y $ is denoted by $ |x-y| $. 
Let $ V $ be an arbitrary point set in $ {\mathbb R}^n $. Let $ a > 0 $ be a real number. By a {\it distance graph} with set of vertices $ V $,
we mean the graph $ G = (V,E) $ whose set of edges $ E $ contains {\it all} pairs of points from $ V $ that are at the distance $ a $ apart:
$$
E = \{\{x,y\}: ~ |x-y|=a\}.
$$

Distance graphs are among the most studied objects of combinatorial geometry. First of all, they are at the ground of the classical Hadwiger--Nelson 
problem, which was proposed around 1950 (see \cite{Had1}, \cite{Soi}) and consists in determining the {\it chromatic number of the space}:
$$
\chi({\mathbb R}^n) = \min \left\{\chi: ~ {\mathbb R}^n = V_1 \sqcup \ldots \sqcup V_{\chi}, ~ \forall ~ i ~~ \forall ~ {\bf x}, {\bf y} \in V_i ~~
|{\bf x}-{\bf y}| \neq 1\right\},
$$
i.e., the minimum number of colors needed to color all the points in $ {\mathbb R}^n $ so that any two points at the 
distance 1 receive different colors. In other words, it is the chromatic number of the unit distance graph whose vertex set coincides with 
$ {\mathbb R}^n $. 

Due to the extreme popularity of the subject, colorings of unit distance graphs are very deeply explored. Let us just refer the reader to several 
books and survey articles \cite{AP,BMP,Ch,KW,Rai1,Rai2,Rai8,Szek}. In particular, the best known lower bounds for the chromatic numbers in dimensions 
$ \le 12 $ are given below:  
$$
\chi({\mathbb R}^2) \ge 4 ~ \cite{Rai1}, ~ \chi({\mathbb R}^3) \ge 6 ~ \cite{Nech}, ~ \chi({\mathbb R}^4) \ge 9 ~ \cite{Exoo}, ~
\chi({\mathbb R}^5) \ge 9 ~ \cite{Cant}, ~ \chi({\mathbb R}^6) \ge 11 ~ \cite{Ci}, ~ \chi({\mathbb R}^7) \ge 15 ~ \cite{Rai1}, ~ 
$$
$$
\chi({\mathbb R}^8) \ge 16 ~ \cite{LR}, ~ \chi({\mathbb R}^9) \ge 21 ~ \cite{KupRai}, ~ \chi({\mathbb R}^{10}) \ge 23 ~ \cite{KupRai}, ~
\chi({\mathbb R}^{11}) \ge 25 ~ \cite{Kup2}, ~ \chi({\mathbb R}^{12}) \ge 27 ~ \cite{Kup1}.
$$
Recently further improvements were announced:
$$
\chi({\mathbb R}^8) \ge 19 ~ \cite{K}, ~ \chi({\mathbb R}^{10}) \ge 26 ~ \cite{Ex}, \cite{K}, ~ \chi({\mathbb R}^{11}) \ge 32 ~ \cite{K}, ~
\chi({\mathbb R}^{12}) \ge 36 ~ \cite{Ex}.
$$
These improvements are essentially based on computer calculations.

In growing dimensions, the following bounds are the best known:
$$
\cite{Rai7} ~~ (1.239\ldots + o(1))^n \le \chi({\mathbb R}^n) \le (3+o(1))^n ~~ \cite{LR}.
$$

In this paper, we consider a special sequence of graphs defined in the following way.

Let $V_n$ be the set of all vectors $v$ from $\mathbb {R}^n$ with coordinates in $\{-1,0,1\}$ and $|v| = \sqrt{3}$. The set $V_n$ can be 
considered as the set of vertices of a graph $G_n = (V_n, E_n)$, where an edge connects two vertices if and only if the corresponding vectors have 
scalar product $1$. Note that $G_1$ and $G_2$ are empty and $G_3$ is just a cube.

Recall that an {\it independent set} in a graph is any set of its vertices which are pairwise non-adjacent and  
the \textit{independence number} of $G$ denoted by $\alpha(G)$ is the size of a maximum independent set in the graph $G$.

\begin{theorem} 
For $n \geq 1$, let $c(n)$ denote the following constant:
\[
c(n) = \begin{cases} 
0    & \mbox{ if } ~ n \equiv 0 \\ 
1    & \mbox{ if } ~ n \equiv 1 \\
2    & \mbox{ if } ~ n \equiv 2 \mbox{ or } 3
\end{cases} \pmod{4}.
\]
Then, the independence number of $G_n$ is given by the formula
$$\alpha(G_n) = \max \{6n-28, 4n - 4c (n)\}.$$
\label{theo1}
\end{theorem}

Actually, the result of Theorem~\ref{theo1} is a far-reaching generalization of a much simpler lemma proved by Zs. Nagy (see \cite{nag}) in 1972 and used 
not only in combinatorial geometry, but also in Ramsey theory. In this lemma, $ G_n' = (V_n', E_n') $, where $ V_n' $ is 
the set of all vectors $v$, $|v| = \sqrt{3}$, with coordinates in $\{0,1\}$ and again an edge connects two vertices if and only if the corresponding vectors have 
scalar product $1$. Lemma states that in this case $ \alpha(G_n') = n - c(n) $.  

The proof of Theorem~\ref{theo1} is given in the following parts: some examples showing the lower bound in Theorem~\ref{theo1} and some preliminaries are given in 
Section \ref{Pre}; the upper bound is proved in Section~\ref{ng14} for $n > 13$ and in Section~\ref{nl14} 
for other values of $n$. Note that, roughly speaking, the quantity 13 is a threshold where the bound $ 6n - 28 $ starts dominating the bound $ 4n $. 


\vskip+0.2cm

As a corollary of Theorem~\ref{theo1} we get the following bounds for the chromatic numbers of Euclidean spaces.

\begin{theorem}
Let $c(n)$ be the constant defined in Theorem~\ref{theo1}. Then, for all $n \geq 3$, we have
$$\chi (\R^n) \geq \chi (\Q^n) \geq \chi(G_n) \ge \frac{|V_n|}{\alpha(G_n)} = \frac{8C_n^3}{\max \{6n - 28 , 4n - c(n)\}}.$$
\label{theo2}
\end{theorem}

Asymptotically, the bound in this theorem is 
$\frac {2}{9}n^2 (1 + o(1))$, which is a weak result. On the other hand, for small values of $n$, the theorem gives the best known bounds, namely:
$$\chi (\mathbb{R}^9) \geq \chi (\mathbb{Q}^9) \geq 21,$$ $$\chi (\mathbb{R}^{10}) \geq \chi (\mathbb{Q}^{10}) \geq 30,$$ $$\chi (\mathbb{R}^{11}) \geq 
\chi (\mathbb{Q}^{11}) \geq 35,$$ $$\chi (\mathbb{R}^{12}) \geq \chi (\mathbb{Q}^{12}) \geq 37.$$

Actually, we will show in Section \ref{prop1} the following stronger result for $n = 9$. 
\begin{proposition}
\label{n9}
The inequalities hold
$$\chi (\mathbb{R}^9) \geq \chi (\mathbb{Q}^9) \geq 22.$$
\end{proposition}

\vskip+0.2cm

\section{Lower bounds in Theorem \ref{theo1} and some preliminaries} \label{Pre}

\subsection{Auxiliary definitions}

Consider the graph $ G_n $. Any of its vertices has three non-zero coordinates and $ n-3 $ coordinates equal to 0. We call {\it base} the set of 
non-zero coordinates of a vertex. To make our exposition more concise, we will use the word ``place'' instead of the word ``coordinate'' or instead the expression
``coordinate position''. For example, it will be convenient to say (a bit informally) ``vertex $v$ intersects place $x$'', if  
the vector $v$ from $\mathbb {R}^n$ corresponding to this vertex has nonzero value of the coordinate $v_x$. For the same reasons, we introduce the 
notion of a {\it signplace}: it is a coordinate with a fixed sign (plus or minus). In particular, from now on, we can say (again, a bit informally)
``vertex $v$ intersects signplace $x^+$ ($x^-$)'', if it has the value of the coordinate $v_x$ equal to $ +1 $ ($ -1 $). Finally, we define the 
{\it degree} of a place (signplace) in a set $ W $ of vertices of $ G_n $ as the number of vertices from $ W $ intersecting this place (signplace).

\subsection{Constructions of independent sets in $ G_n $}

It suffices to show that $ \alpha(G_n) \ge 4n - 4c(n) $ and that $ \alpha(G_n) \ge 6n-28 $. 

The first construction is as follows. Consider the first 4 places. Take all the 4 bases that can be taken on these places. For each base, consider 4 variants:
$$ 
1,1,1; ~~~ 1,-1,-1; ~~~ -1,1,-1; ~~~ -1,-1,1.
$$
Clearly any two vectors with these bases have scalar product different from 1. We call this construction (and its natural analogs) {\it quad}. 

Take $ [n/4] $ consecutive quads. If the remainder still consists of 3 places, then add 4 more bases. Eventually, we get exactly $ 4n-4c(n) $ vectors that form 
an independent set in $ G_n $. 

Now, let us make the second construction. Take the following vectors:
$$
1,-1,0,1,0, \dots, 0,0,0,0; ~~ 1,-1,0,0,1,0, \dots, 0,0,0,0; ~~ \dots; ~~  1,-1,0,0,0, \dots, 0,1,0,0,0;
$$  
$$
0,1,-1,1,0, \dots, 0,0,0,0; ~~ 0,1,-1,0,1,0, \dots, 0,0,0,0; ~~ \dots; ~~ 0,1,-1,0,0, \dots, 0,1,0,0,0;
$$  
$$
-1,0,1,1,0, \dots, 0,0,0,0; ~~ -1,0,1,0,1,0, \dots, 0,0,0,0; ~~ \dots; ~~  -1,0,1,0,0, \dots, 0,1,0,0,0.
$$  
In each line, we have a set of vectors, which is a particular case of what we will call {\it snake} in Section \ref{ng14} and later. In every snake, we have 
$ n-6 $ vertices. Thus, the total amount of vertices here is $ 3n - 18 $. Obviously, the union of these snakes is an independent set in $ G_n $. Moreover, 
we can add to it 4 more vectors, which have a common base --- the three first places: say, 
$$ 
1,1,1,0,\dots,0; ~~~ 1,-1,-1,0,\dots,0; ~~~ -1,1,-1,0,\dots,0; ~~~ -1,-1,1,0,\dots,0.
$$
The whole construction is a particular case of a {\it cobra} discussed later in more details. Here the cobra contains $ 3n-14 $ vertices.

Of course, we can take one more cobra, whose ``head'' is on the three last places and whose ``tail'' consists of minus ones instead of ones. Eventually, we 
get exactly $ 6n-28 $ vertices forming an independent set in $ G_n $.  

The lower bound is proven. 

\vskip+0.3cm

It is worth noting that in Section \ref{ng14} we will make a rather subtle analysis of possible independent sets in $ G_n $. One would be able to derive from 
this analysis a complete description of examples giving the lower bound in Theorem \ref{theo1}. However, we will not present such description explicitly in this 
paper. 

It is also worth noting that in the above example having $ 6n - 28 $ vertices and avoiding the scalar product 1, the scalar product $ - 3 $ is also absent. Moreover, one 
can exclude 6 vertices from that example so that the scalar product $-2$ disappears as well. 

\subsection{Basic lemma}

Let $ A $ be an arbitrary independent set of the maximum size in $ G_n $. We already know that $|A| \geq  \max \{6n-28, 4n - 4c (n)\}$. 
Assume that we exclude some signplaces and all the vertices from the graph $ G_n $ intersecting them. Then we get a new graph $ G' $ with a possibly smaller 
independent set $ A' $. Denote by $ a(A') $ the maximum degree of a signplace in the set $ A' $. Denote by $ m(A') $ the number of signplaces in $ A' $. 

The following lemma is an important ingredient in the proof of the upper bound. 

\begin{lemma}  
Assume that we exclude $k$ signplaces. Assume that the number of vertices excluded from $ A $ does not exceed $2k$. Then we have either $a(A') \geq 5$ or $m(A') < 14$. 
\label{lem1}
\end{lemma}

\paragraph{Proof of the lemma.} By pigeon-hole principle $a(A') \geq 3(|A| - 2k)/(2n - k)$. If $ |A| \ge 4n $, then $ a(A') \ge 6 $ and we are done. 
The inequality $ |A| \ge 4n $ is true for $ n = 8, 12 $ and $ n \ge 14 $. Thus, it remains to consider only $n = 7, 9, 10, 11, 13$. If $ n = 7 $, then 
$ |A| \ge \max \{14,20\} = 20 $. If $ k = 1 $, then $ 3(|A| - 2k)/(2n - k) \ge 54/13 $, i.e., $ a(A') \ge 5 $. If $ k \ge 2 $, then it may happen that 
$ 3(|A| - 2k)/(2n - k) \le 4 $. But in this case, $ m(A') = 2n-k \le 14 - 2 < 14 $. The same argument works for the 4 other values of $ n $. The proof is 
complete.

\section{Proof of Theorem \ref{theo1} in the case $n \geq 14$}\label{ng14}

\subsection{Starting the proof}

Let $ A $ be an arbitrary independent set of the maximum size in $ G_n $. Assume that we have already excluded several signplaces with the 
corresponding vertices (see Section 2.3). By Lemma \ref{lem1} either $a(A') \geq 5$ or $m(A') \leq 13$. The second case will be considered in Section 3.4.
So we assume that $a(A') \geq 5$.  

Consider a signplace with the maximum degree. Call it $x_1^+$ (each time when we choose a sign we can choose plus without loss of generality) and consider the set of 
vertices intersecting it (denote it by $N_{x_1^+}$). Note that no base can contain more than two vertices from $N_{x_1^+}$. Thus, we have at least three 
different bases. Also it is clear that any two bases containing vertices from $N_{x_1^+}$ intersect in exactly two signplaces. There are two different possibilities.

\begin{enumerate}
\item Among the bases, we have $\{x_1, x_2, x_3\}, \{x_1, x_2, x_4\}, \{x_1, x_3, x_4\}$.
This case will be referred to as ``quad'' (cf. Section 2.2). \label{case1}
\item All the bases contain both $x_1$ and $x_2$.
This case will be referred to as ``snake'' (cf. Section 2.2). \label{case2}
\end{enumerate}

The formal definition of a quad will be given in the next section, where we will analize Case \ref{case1}. The same is for a snake in Section 3.3. In Section 3.4,
we will complete the proof.

\subsection{The first case --- ``quad''}

We know that $ a(A') \ge 5 $. At the same time, $a(A') \leq 6$, since otherwise the vertices from $N_{x_1^+}$ use at least 4 bases and therefore there is a base among
$\{x_1, x_2, x_3\}, \{x_1, x_2, x_4\}, \{x_1, x_3, x_4\}$ such that it intersects the fourth base only on $x^+_1$, which is impossible. Put $ a = a(A') $. 

Thus, we have exactly three bases containing the vertices from $N_{x_1^+}$. Two of them (without loss of generality $\{x_1, x_2, x_3\}, \{x_1, x_2, x_4\}$)
contain exactly $2$ vertices from $N_{x_1^+}$ each, and the third one contains at least $1$ vertex. Since 
$ \{x_1^+, x_2, x_3\} $ contains two vertices, it intersects all the four signplaces in $ x_2, x_3 $; the same holds for 
$ \{x_1^+, x_2, x_4\} $, which means that all the six signplaces of $ x_2, x_3, x_4 $ are necessarily intersected. 

Consider the set $ U $ of 
all vertices intersecting $\{x_1, x_2, x_3, x_4 \}$. There could be the following possibilities.

\begin{itemize}

\item Some vertices from $ U $ intersect $x_1$. There are at most $2a$ such vertices.

\item Some vertices from $ U $ lie on the base $\{x_2, x_3, x_4\}$. There are at most $4$ such vertices.

\item Some vertices from $ U $ intersect $\{x_2, x_3, x_4\}$ in one place and are not counted above. 
Actually, there are no such vertices because for every signplace in $\{x_2, x_3, x_4\}$ 
a vertex with a base in $\{ x_1, x_2, x_3, x_4 \}$ exists (do not forget that $\{x_1, x_2, x_3\}, \{x_1, x_2, x_4\}$
contain exactly $2$ vertices each, and the third base contains at least $1$ vertex).

\item Some vertices from $ U $ intersect $\{x_2, x_3, x_4\}$ in two places and are not counted above. Again, there are no such vertices. Indeed,
assume that some vertex (call it $v$) intersects $\{x_2, x_3, x_4 \}$ 
in $\{ x_i, x_j \}$. Then  $\{x_1, x_2, x_3\}$ or $\{x_1, x_2, x_4\}$ intersects $\{ x_i, x_j \}$ in exactly one place. This is impossible, since
we know that two vertices from $N_{x_1^+}$ lie on $\{x_1, x_2, x_3\}$ and two vertices from $N_{x_1^+}$ lie on $\{x_1, x_2, x_4\}$.

\end{itemize}

Summarizing, we have at most $2a + 4 \leq 16$ vertices intersecting $8$ signplaces. We call any of the corresponding constructions {\it quad}.

\vskip+0.2cm

Now we may assume that $ A $ was transformed into $ A' $ in the following way (more details will be given in Section 3.4).

\begin{itemize}

\item First, all the signplaces of degree less than $3$ have been deleted one by one. Note that by Lemma~\ref{lem1} during this process either $a \geq 5$ or $m \leq 13$.

\item Second, all the quads have been deleted one by one. Note that again by Lemma~\ref{lem1} during this process either $a \geq 5$ or $m \leq 13$ (everytime
the number of excluded signplaces is 8 and the number of excluded vertices is at most 16). 

\item Third, once again, all the signplaces of degree less than $3$ have been deleted one by one. 
Obviously, there are no new quads and still by Lemma~\ref{lem1} $a \geq 5$ or $m \leq 13$.

\end{itemize}

As before, we assume that $ a \geq 5 $ (since the case $ m \le 13 $ is considered in Section 3.4), and 
so we are prepared to the next case, in which we have $a(A') \geq 5$, there are no quads, and every signplace has degree at least $3$.

\subsection{The second case --- ``snake''}

We start with a formal definition of a snake.

\begin{definition}  \textbf{Snake} is a set of vertices intersecting a signplace and a place and containing at least $5$ vertices. \textbf{Head} of a snake is a 
couple of places, which intersect every vertex, and \textbf{tail} of a snake is the set of the remaining signplaces in each vertex. 
\textbf{Size} of a snake is the number of its vertices.
\end{definition}

Clearly in the current case we have a snake of size $a \geq 5$ in $ A' $. 
Let it be based on $\{ x_1^+, x_2 \}$ (with the head being $\{ x_1, x_2 \}$). 
Note that the size of its tail is equal to $a$, since vertices cannot intersect on tail. 

Our aim is to prove that we can exclude some $ t $ signplaces with at most $ 3t-14 $ vertices. Moreover, we will show that there is a special construction 
(``cobra'', cf. Section 2.2), which has exactly $ 3t-14 $ vertices on $ t $ signplaces and which is the only such construction up to the graph symmetries. 

We have an alternative.

\begin{enumerate}

\item  We can exclude $4 + a$ signplaces ($\{ x_1, x_2 \}$ with all possible signs and $a$ signplaces of the tail) and $3a-2$ vertices. 

\item We have at least $3a-1$ vertices intersecting the signplaces mentioned in the previous point. 

\end{enumerate}

In the first case, our aim is realized, since we can put $ t = 4+a $ and get $ 3t-14 = 3a-2 $. In the second case, the analysis will be much longer.

Let us consider the second case of the alternative. Each vertex intersecting the tail of the snake that we analize should 
intersect the head as well, and each of the $a$ initial vertices intersects the head on two signplaces. 
Hence the sum of the degrees of the head signplaces is at least $4a-1$. But there is no signplace with degree exceeding $a$, 
so the degrees of the signplaces in the head are either $$a, a, a, a ~~~ \text{or} ~~~ a, a, a, a-1.$$ Anyway we have a place with two signplaces of degree exactly $a$. 
Without loss of generality, this place is $x_1$. Since all quads are already excluded, we have two snakes with signplaces on $x_1$: 
one signplace is $x_1^+$ and the second one is $x_1^-$. Consider their heads. They could both lie on $\{x_1, x_2\}$, or they could lie on $\{x_1, x_2\}$ 
and $\{x_1, x_3\}$ respectively.

In the first case, all the four signplaces of the head have degree $a$ solely due to $2a$ vertices from the snakes. 
In addition, there are vertices intersecting the tail (since 
the degree of each signplace is at least 3 and two snakes could provide only two vertices on a signplace). Each vertex intersecting the tail should intersect the 
head as well, so the degree of some signplace in the head exceeds $a$, which contradicts the assumption that $a$ is the maximum value of the degree.

We are left with the second case: there are two snakes of size $a$ with heads on $\{ x_1, x_2\}$ and $\{ x_1, x_3 \}$.

Let $Q$ be the set of vertices lying fully on base $\{ x_1, x_2, x_3 \}$. Denote by $B$ the set of signplaces in the intersection of the tails. Let $C_1$ and $C_2$ 
be the sets of the remaining signplaces in the corresponding tails. Let $q, b, c_1, c_2$ be the sizes of the corresponding sets.
We have already described all the vertices intersecting $x_1$, since the maximum degree is equal to $a$. Consider the sum of the degrees of 
the signplaces on $x_2$ and $x_3$. 
Since the degree of each signplace is at least 3, we have a new vertex for each signplace from the intersection of the tails. 
Each vertex of this type should intersect both heads, and it cannot contain $x_1$. Therefore, it contains both $x_2$ and $x_3$ and adds $2$ to our sum. We have at least 
two vertices intersecting each signplace of the symmetric difference of the tails. Each vertex of this type should intersect the head of a corresponding snake and could 
intersect two signplaces of its tail. In total, these vertices add at least $$2(c_1 + c_2)/2 = c_1 + c_2$$ to the sum. Each of the $2a$ initial vertices intersects 
$\{x_2, x_3\}$. Each vertex from $Q$ adds yet another $1$ to the sum, since it intersects $\{x_2, x_3\}$ on two places.
Again in total, the sum of the degrees of the four signplaces on $\{x_2, x_3 \}$ is at least $$2b + c_1 + c_2 + 2a + q.$$ On the other hand, 
since the degree of each signplace 
is at most $a$, this sum does not exceed $4a$. So we have $$2b + c_1 + c_2 + q \leq 2a.$$ Each vertex from $Q$ is in one snake. Consequently, 
$q = q_1 + q_2$ ($q_i$ is the 
number of vertices lying in a corresponging snake), $$b + q_1 + c_1 = a, ~~~ b + q_2 + c_2 = a,$$ and the inequality always turns to equality! 

Thus, there is a set of $$t := 6 + b + c_1 + c_2$$ signplaces intersected by $$q + 3b + 3c_1 + 3c_2 \leq 3(b+c_1+c_2) + 4 = 3t - 14$$ vertices. 

The second case of the alternative is complete, and our aim is attained. However, we will also prove below an upper bound on $t$. 

Suppose that the number of vertices is exactly $3t-14$. It means that all the intermediate inequalities turned to equalities. 
The last inequality turns to equality only when $q = 4$. One can see that any vertex intersecting a signplace
from $C_1$ or $C_2$ should intersect $2$ vertices of the tail, so it intersects $\{x_1, x_2, x_3\}$ only on $1$ vertex, 
which contradicts $q = 4$. Hence $C_1 = C_2 = \emptyset$ and $b = a - 2$. 
For every signplace from $x^+ \in B$ there is a vertex intersecting $x^+$ and lying on base $\{x, x_2, x_3\}$. It turns out that there is a third snake on $\{x_2, x_3\}$.  
We call {\it cobra} the union of such 3 snakes. Finally, one can see that there is no place $x$ such that $x^+, x^-$ lie in $B$, otherwise
there is an edge between two vertices on the base $\{x_1, x_2, x\}$.

As a result, $t$ does not exceed $n+3$, and the tail of a corresponding snake can not contain two signplaces on the same place.

Summing up the above, if we have no quad, then there is a cobra, which consists of three snakes with a common tail and pairwise intersecting heads.
It has $3t - 14$ vertices on $t$ signplaces, $8 \leq t \leq n+3$.

\subsection{Finishing the proof} \label{5.4}

In the previous sections we have shown that there are the following options.

\begin{itemize}
\item To exclude a signplace and at most $2$ vertices intersecting it.
\item To exclude $8$ signplaces and at most $16$ vertices intersecting it.
\item To exclude $t$ signplaces with at most $3t - 14$ vertices ($8 \leq t \leq n+3$).
\item To get $ m(A') \le 13 $. 

\end{itemize}

Clearly the two first options yield at most $2m$ vertices on $m$ signplaces. The same is true for the fourth option: we show it in Appendix using 
computer simulations. 

According to this, only the following cases could occur.

\begin{enumerate}

\item There is no cobra. Then the number of vertices does not exceed $4n \leq 6n-28$.

\item There is one cobra and $t \leq n$.
Then the number of vertices does not exceed 
$$3t - 14 + 2(2n - t) \leq 5n - 14 \leq 6n - 28.$$

\item There is one cobra and $t = n+1$.
We are left to prove that $n-1$ signplaces on $n-3$ places can contain at most $2n - 3$ vertices.
Suppose the contrary, then, by pigeon-hole principle, there is a signplace of degree at least $3(2n-2)/(n-1) = 6$. Using the same arguments as 
in Sections 3.1--3.3 we get a quad or a cobra, but both constructions contain $3$ places with $2$ signplaces, which contradicts our assumptions.

\item There is one cobra and $t = n+2$.
We are left to prove that $n-2$ signplaces on $n-3$ places can contain at most $2n - 6$ vertices.
Then the number of vertices does not exceed $$3(n+2) - 14 + 2n - 6 \leq 5n - 14\leq 6n - 28.$$
Again, suppose the contrary, so there is a singplace of degree at least $3(2n-5)/(n-2)$. For $n < 4$ the claim is obvious, and
for $n \geq 4$ we have $3(2n-5)/(n-2) > 4$. Using the same arguments as in Sections 3.1--3.3
we get a quad or a cobra, but both constructions contain $3$ places with $2$ signplaces. Thus, we get a contradiction.

\item There is one cobra and $t = n+3$.
All other signplaces lie on distinct places, and we can apply Nagy's lemma (see~Section 1 and \cite{nag}) to get an upper bound $n-3$ for the number of vertices.
Then the number of vertices does not exceed $$3(n+3) - 14 + n - 3 = 4n - 8 < 5n - 14 \leq 6n - 28.$$ 

\item There are two or more cobras. Then the bound $6n - 28$ for the number of vertices is straightforward. 

\end{enumerate}

The proof of Theorem \ref{theo1} for $ n \ge 14 $ is complete.

\section{Proof of Theorem 1 in the case $n < 14$}\label{nl14}

We will deal with small values of $n$ one by one.  We have already proven the lower bounds in Section 2.2. So we have to prove only the upper bounds.
Let us start from simple cases. 

\begin{itemize}
\item In the case $n \leq 6$ the answer is easily found on a computer.
\item The case $n = 7$ was considered by Cibulka in \cite{Ci}. 
\item In the case $n = 9$ we have to show that $\alpha (G_9) \leq 32$, and this will follow from the case $n = 10$ below, since $\alpha (G_n)$ is non-decreasing.
\end{itemize}

In the remaining cases the main lines of the proofs will be the same as in Sections 3.1--3.3, but we will make appropriate changes in the 
ends of the proofs.

\subsection {The cases $n = 8$ and $n=12$}

We state that $\alpha (G_8) \leq 32$ and $\alpha (G_{12}) \leq 48$.
Suppose the contrary. It means that $\alpha (G_n) > 4n$ for some $n \in \{8, 12\}$. Hence, by Section~\ref{5.4} 
there is a cobra on $t$ signplaces with $t \geq 15$, since otherwise $3t-14 \leq 2t$ and thus in all the possibilities from Section~\ref{5.4} we have at
most $2m$ vertices on $m$ signplaces. But a cobra has at most $n+3$ signplaces, so the only option is $n = 12$, $t = 15$. 
Let us exclude the cobra and note that all the $9$ remaining signplaces lie on different places. 
Now we can apply Nagy's lemma (see~Section 1 and \cite{nag}) to get in this case an upper bound $3t - 14 + 8 = 39 < 48$.
We get a contradiction and the claim is proven.

\subsection{The case $n = 10$}

We are going to prove that $\alpha (G_{10}) \leq 32$. Suppose we have an independent set $A$ whose size is at least $32$.
Let us follow the proof of the case $n \geq 14$. If there is a quad, we exclude it and refer to the Appendix, where
the results of some computer calculations are given, in particular, for $m = 12$, $ l = 6 $ ($ m $ is 
the number of signplaces, and $ l $ is the number of places). Thus, we get the needed bound by $ 16 + 16 = 32 $.

Then, we exclude some $k$ signplaces whose degrees are at most 2. Suppose that the maximum degree in the remaining graph is at most $4$.
By Lemma~\ref{lem1} we get the bound $k \geq 8$. Then $|A| \leq 2k + 4 (2n - k)/3$. Hence $k \geq 10$. If among the
excluded signplaces there were $6$ signplaces on $3$ places, then we can apply the bound for $G_7$: $|A| \leq \alpha (G_7) + 2 \cdot 6 = 32$. 
Otherwise there are at least $n-2 = 8$ different places (and thus $k$ is at most $12$) in the remaining graph, and
so one can apply computer calculations (see Appendix) for $m = 2n - k \leq 10$, $l \geq 8$ and get the bounds by $2\cdot 10 + 12 = 32$ or $2\cdot 11 + 10 = 32$ 
when $k = 10$ or $11$, respectively. For the case $k = 12$ Nagy's lemma gives us the inequality $|A| \leq 2\cdot 12 + 8 = 32$.

If a cobra exists, it cannot be intersected by any vertex (including the already deleted ones). The cobra uses at least $8$ signplaces and 
leaves at least $n-3 = 7$ places, which are not fully used (i.e., at least one of the signplaces on a place is free). Let $t$ be the number of 
signplaces occupied by the cobra and $m$ be equal to $2n - t$. Then a graph $G'_{10}$, which is obtained by excluding the signplaces 
intersected by the cobra, has $m$ signplaces lying on $l \geq 7$ different places.

Now we only have to sum up the number of vertices from the cobra and the 
number of vertices on the remaining signplaces (the last number is given in Appendix). So let us consider all the cases with $t$ varying from $8$ to $n+3 = 13$:

\begin{itemize}

\item $t = 8,  m = 12, 3t - 14 + 22 = 32$;
\item $t = 9,  m = 11, 3t - 14 + 19 = 32$;
\item $t = 10, m = 10, 3t - 14 + 16 = 32$;
\item $t = 11, m = 9,  3t - 14 + 13 = 32$, since we have at least $7$ different places;
\item $t = 12, m = 8,  3t - 14 + 10 = 32$, since we have at least $7$ different places;   
\item $t = 13, m = 7,  3t - 14 + 5  = 30$, since all the places are distinct.

\end{itemize}

\subsection{The case $n = 11$}

We are going to prove that $\alpha (G_{11}) \leq 38$. Assume the contrary.
If we have a quad, then after its exclusion we get exactly $G_7$, but $\alpha (G_7) = 20$ and $20 + 16 < 38$.

Now exclude one by one all the vertices with degree at most 2. Denote the number of excluded vertices by $k$. Following the proof of Lemma~\ref{lem1} we see that if 
$k < 13$, then $a > 4$, which means that a cobra exists. First, let $k \geq 13$.
If among the excluded signplaces there where $8$ signplaces on $4$ places, then we can once again apply the bound for $G_7$. Otherwise there are at
least $n-3 = 8$ different places in the remaining graph, so we can refer to Appendix with $m = 2n-k \le 9$, $l \geq 8$ and get
at most $38$ vertices in $A$.

If a cobra exists, it cannot be intersected by any vertex (including the already deleted ones). The cobra uses at least $8$ signplaces and 
leaves at least $n - 3 = 8$ places, which are not fully used. Let $t$ be the number of 
signplaces occupied by the cobra and $m$ be equal to $2n - t$. Then a graph $G'_{11}$, which is obtained by excluding the signplaces 
intersected by the cobra, has $m$ signplaces lying on $l \geq 8$ different places.

If $t \leq 8$, then we use the fact that
all the signplaces, which are not in a cobra, have degrees at most two and get the declared bound: $3t - 14 + 2 (22 - t) = t + 30 \leq 38$.

Now it remains to consider the cases with $t$ from $9$ to $n+3 = 14$ and use the results from Appendix:

\begin{itemize}
\item $t = 9,  m = 13, 3t - 14 + 25 = 38$;
\item $t = 10, m = 12, 3t - 14 + 22 = 38$;
\item $t = 11, m = 11, 3t - 14 + 19 = 38$;
\item $t = 12, m = 10, 3t - 14 + 16 = 38$;  
\item $t = 13, m = 9,  3t - 14 + 13 = 38$, since we have at least $8$ different places;
\item $t = 14, m = 8,  3t - 14 + 10 = 38$, since all the places are distinct. 
\end{itemize}

\subsection{The case $n = 13$}

We are going to prove that $\alpha (G_{13}) \le 50$. Assume the contrary.
If we have a quad, then after its exclusion we get exactly $G_9$, but $\alpha (G_9) = 32$ and $32 + 16 < 50$.

Now exclude one by one all the vertices with degree at most 2. Denote the number of excluded vertices by $k$. Following the proof of Lemma~\ref{lem1} we see that
if $k < 23$, then $a > 4$, which means that a cobra exists. The case $k \geq 23$ is obvious.

Thus, a cobra exists, and it cannot be intersected by any vertex (including the already deleted ones). Define $ t $ as in the pevious cases. If $t < 13$, then 
we use the fact that all the signplaces, which are not in a cobra, have degrees at most two and get the declared bound: $3t - 14 + 2 (26-t) = t + 38 < 50$.

Now it remains to consider the cases with $t$ from $13$ to $n+3 = 16$ and use the results from Appendix:

\begin{itemize}

\item $t = 13, m = 13,  3t - 14 + 25 = 50$; 
\item $t = 14, m = 12,  3t - 14 + 22 = 50$; 
\item $t = 15, m = 11,  3t - 14 + 19 = 50$;  
\item $t = 16, m = 10,  3t - 14 + 16 = 50.$

\end{itemize}

The proof is complete.

\section{Proof of Proposition~\ref{n9}}\label{prop1}

Suppose that $\chi (G_9) = 21$. Clearly $\frac {|V(G_9)|}{\alpha (G_9)} = 21$, and therefore every independent set has size $32$.
Revising the proof of Theorem~\ref{theo1} we see that the only way to reach $32$ vertices in an independent set is by taking a couple of full quads.
Thus, we have a collection of $21$ pairs of full quads (denote it by $A$); this collection covers each base exactly two times, since every full quad has exactly 
$4$ vertices on every covered base.
Note that every pair of quads does not cover exactly one place, so one can split $A$ into nine disjoint parts: $$A = A_1 \sqcup \ldots \sqcup A_9.$$
Let $S_1$ be the set of all bases such that each of them does not contain the first place. Obviously $|S_1| = C^3_8 = 56$. Consider a pair of quads $p \in A$. 
Note that $p$ covers $8$ bases from $S_1$, if $p \in A_1$, and $5$ bases from $S_1$ otherwise. Denote the cardinalities of $A_1$ and $A \setminus A_1$ by $a$ and 
$b$ respectively.
Every set in $S_1$ is covered twice, and therefore we have $2|S_1| = 112 = 8a + 5b$. Hence there are the following possibilities:  $(a = 14, b = 0)$,  $(a = 9, b = 8)$ 
and $(a = 4, b = 16)$. But $a + b = |A| = 21$, so we get a contradiction.

Proposition~\ref{n9} is proved.

\section*{Appendix. Computer calculations}\label{comp}

Let $F$ be a subgraph of $G_n$. Denote the number of signplaces and places intersecting the vertices of $F$ by $m$ and $l$ respectively.

We use the standard Bron--Kerbosch algorithm (see \cite{BK}) and get the following results:

\begin{itemize}

\item if $m = 13$, then $\alpha (F) \leq 25$;
\item if $m = 13$ and $l = 7$,  then $\alpha (F) \leq 18$;
\item if $m = 12$, then $\alpha (F) \leq 22$;
\item if $m = 12$ and $l = 6$, then $\alpha (F) \leq 16$;
\item if $m = 11$, then $\alpha (F) \leq 19$;
\item if $m = 10$, then $\alpha (F) \leq 16$;
\item if $m = 10$ and $l \geq 8$, then $\alpha (F) \leq 12$;
\item if $m = 9$,  then $\alpha (F) \leq 16$;
\item if $m = 9$ and $l \geq 7$,  then $\alpha (F) \leq 13$;
\item if $m = 9$ and $l \geq 8$,  then $\alpha (F) \leq 10$;
\item if $m = 8$,  then $\alpha (F) \leq 16$;
\item if $m = 8$ and $l \geq 7$,  then $\alpha (F) \leq 10$;
\item if $m = 7$,  then $\alpha (F) \leq 10$;
\item if $m = l = 7$,  then $\alpha (F) \leq 5$;
\item if $m = 6$,  then $\alpha (F) \leq 7$;
\item if $m = 5$,  then $\alpha (F) \leq 5$;
\item if $m = 4$,  then $\alpha (F) \leq 4$;
\item if $m = 3$,  then $\alpha (F) \leq 1.$

\end{itemize}

\vskip+1cm

\noindent\textbf{Acknowledgements.} 
Danila Cherkashin was supported by the grant of the Russian Scientific Foundation (grant N16-11-10039). 
All coauthors are grateful to Vyacheslav Sokolov for computer calculations.

\end{document}